\documentclass[12pt]{amsart}

\usepackage{epsfig,amssymb}

\newcommand{\be}{\begin{equation}} 
\newcommand{\ee}{\end{equation}}
\newcommand{\beq}{\begin{eqnarray}}
\newcommand{\eeq}{\end{eqnarray}}
\newcommand{\bt}{\beta}
\newcommand{\bl}{\begin{lemma}}
\newcommand{\el}{\end{lemma}}
\newcommand{\bm}{\begin{pmatrix}}
\renewcommand{\em}{\end{pmatrix}}
\newcommand{\bml}{\begin{multline}}
\newcommand{\eml}{\end{multline}}
\newcommand{\ba}{\begin{array}}
\newcommand{\ea}{\end{array}}

\newcommand{\la}{\label}
\newcommand{\ci}{\cite}

\newcommand{\de}{\delta}

\newcommand{\al}{\alpha}
\newcommand{\ga}{\gamma}

\newcommand{\si}{\sigma}

\newcommand{\om}{\omega}

\newcommand{\lb}{\lambda}

\renewcommand{\th}{\theta}
\newcommand{\Th}{\Theta}

\newcommand{\ep}{\varepsilon }

\newcommand{\bi}{\bibitem}

\newfont{\msbm}{msbm10 scaled\magstep1}
\newfont{\msbms}{msbm7 scaled\magstep1} 

\newcommand{\bbr}{\mbox{$\mbox{\msbm R}$}}

\newcommand{\bbc}{\mbox{$\mbox{\msbm C}$}}

\newcommand{\bbz}{\mbox{$\mbox{\msbm Z}$}}

\newcommand{\esssup}{\mbox{ess sup }}
\newcommand{\essinf}{\mbox{ess inf }}

\newtheorem{theorem}{Theorem}[section]
\newtheorem{lemma}[theorem]{Lemma}

\newtheorem{corollary}[theorem]{Corollary}

\theoremstyle{definition}

\theoremstyle{remark}
\newtheorem{remark}[theorem]{Remark}

\numberwithin{equation}{section}

\begin{document}
\def\wt{\widetilde}
\def\wh{\widehat}
\title[Toeplitz eigenvalues]{Eigenvalues of Toeplitz matrices in the bulk of the spectrum}
\author{P. Deift}
\address{Courant Institute of Mathematical Sciences, New York, NY, USA}
\author{A. Its}
\address{
Indiana University -- Purdue University  Indianapolis,
Indianapolis, IN, USA}
%
\author{I. Krasovsky}
\address{
Imperial College, London, United Kingdom}
%


\dedicatory{To Raghu Varadhan: in admiration for all your achievements}

\begin{abstract}
The authors use results from \ci{DIKt1,DIKt2} to analyze the asymptotics 
of eigenvalues of Toeplitz matrices with certain continuous and discontinuous symbols. 
In particular, the authors prove a conjecture of Levitin and Shargorodsky on the 
near-periodicity of Toeplitz eigenvalues.
\end{abstract}

\maketitle

\section{Introduction}
Let $f(z)$ be a complex-valued function integrable over the unit circle $C$
with Fourier coefficients
\[
f_j={1\over 2\pi}\int_0^{2\pi}f(e^{i\theta})e^{-i j\theta}d\theta,\qquad j=0,\pm1,\pm2,\dots
\]
We are interested in the eigenvalues of $n$-dimensional Toeplitz matrices with symbol $f(z)$,
\be\la{TD}
T_n(f)=(f_{j-k})_{j,k=0}^{n-1}.
\ee
Denote the corresponding Toeplitz determinants
\be
D_n(f)=\det T_n(f).
\ee

Over the years the eigenvalues $\lb^{(n)}_k$, $k=1,\dots,n$, of $T_n(f)$ and
their asymptotics as $n\to\infty$, have been analyzed extensively (for an
outline of the work in this direction, see, for example \ci{BG}, pp 256--259). 
In this paper we discuss the eigenvalues of $T_n(f)$ for large $n$ using general theorems
on Toeplitz determinants proved in \cite{DIKt1,DIKt2}. In principle, we can address
the case where the symbol $f(e^{i\theta})$ has a fixed number of 
Fisher-Hartwig singularities \ci{FH,L}, i.e., 
when $f(e^{i\theta})$ has the following form on the unit circle $C$:
\be\la{fFH}
e^{V(z)} z^{\sum_{j=0}^m \bt_j} 
\prod_{j=0}^m  |z-z_j|^{2\al_j}g_{z_j,\bt_j}(z)z_j^{-\bt_j},\qquad z=e^{i\th},\qquad
\theta\in[0,2\pi),
\ee
for some $m=0,1,\dots$,
where
\begin{eqnarray}
&z_j=e^{i\th_j},\quad j=0,\dots,m,\qquad
0=\th_0<\th_1<\cdots<\th_m<2\pi;&\la{z}\\
&g_{z_j,\bt_j}(z)\equiv g_{\bt_j}(z)=
\begin{cases}
e^{i\pi\bt_j}& 0\le\arg z<\th_j\cr
e^{-i\pi\bt_j}& \th_j\le\arg z<2\pi
\end{cases},&\la{g}\\
&\Re\al_j>-1/2,\quad \bt_j\in\bbc,\quad j=0,\dots,m,&
\end{eqnarray}
and $V(e^{i\theta})$ is a sufficiently smooth function, e.g. $C^{\infty}$, on the unit circle.
We  assume that $z_j$, $j=1,\dots,m$, are genuine singular points, i.e.,
either $\al_j\neq 0$ or $\bt_j\neq 0$. However, we always include $z_0=1$ explicitly 
in (\ref{fFH}), even when $\al_0=\bt_0=0$. The $\bt_j$'s are not uniquely determined 
by the symbol $f(e^{i\th})=f(e^{i\th},\{\al_k\},\{\bt_j\})$. Indeed, if we replace $\{\bt_j\}$
with $\{\wt\bt_j=\bt_j+n_j: n_j\in\bbz, \sum_{j=0}^m n_j=0\}$, then we obtain the relation
\[
 f(e^{i\th},\{\al_k\},\{\bt_j\})=\left(\prod_{k=0}^m z_j^{n_j}\right)
f(e^{i\th},\{\al_k\},\{\wt\bt_j\})
\]
for the symbol $f$. The function $f(e^{i\th},\{\al_k\},\{\wt\bt_j\})$ is called the FH-representation
of the symbol corresponding to $(n_j)_{j=0}^m$.

Our analysis is based on the following observations. First, it is obvious that
the characteristic polynomial of a Toeplitz matrix is a Toeplitz determinant
with symbol shifted by a constant:
\be
\det(T_n(f)-\lb I)=D_n(f-\lb).
\ee
Thus the eigenvalue problem for Toeplitz matrices is equivalent to the problem of looking for the 
zeros of Toeplitz determinants. Let us denote
\be
f(z;0)=f(z),\qquad f(z;\lb)=f(z)-\lb.
\ee
The second observation is that if $f(z)$ is of type (\ref{fFH}) then $f(z;\lb)=f(z)-\lb$ is also
of type (\ref{fFH}) with changed or added Fisher-Hartwig singularities.
For example, as we will see below, if 
\be\label{f1}
f(z)=e^{V(z)},
\ee
where $V(z)$ is a real-valued 
$C^{\infty}$ function on the unit circle such
that $f(e^{i\th})$ is strictly increasing on $\th\in(0,\wt\th)$ for some $\wt\th$,
and strictly decreasing on $(\wt\th,2\pi)$, then for $\min_{z\in C} f(z) 
<\lb<\max_{z\in C} f(z)$, the function
$f(z)-\lb$ will have 2 Fisher-Hartwig singularities with parameters $\al_1=\al_2=1/2$,
$\bt_1=-\bt_2=1/2$ at the points $z_1$, $z_2$ where $f(z)-\lb$ vanishes. 

From these observations, we see that the problem of computing the asymptotic behavior of 
the eigenvalues of $T_n(f)$ reduces to an asymptotic problem for Toeplitz determinants with
Fisher-Hartwig singularities of the kind considered in \ci{DIKt1}, \cite{DIKt2}. 
To leading order, the behavior of the eigenvalues of $T_n(f)$ is determined by the condition
of vanishing of the leading term in the asymptotics of $D_n(f-\lb)$.  
 
An important role in the asymptotic analysis of Toeplitz determinants with symbol (\ref{fFH})
is played by the seminorm
\be\la{norm}
|||\bt|||=\max_{j,k}|\Re\bt_j-\Re\bt_k|,
\ee
where the indices $j,k=0$ are omitted if $z=1$ is not a singular point, i.e. if $\al_0=\bt_0=0$.
If $m=0$, set $|||\bt|||=0$.
If $|||\bt|||<1$, the asymptotics of $D_n(f)$ are given by an explicit formula
(see (1.12) in \cite{DIKt1}) whose leading term is nonzero provided $\al_j\pm\bt_j\neq
-1,-2,\dots$. However, if $|||\bt|||\ge 1$, there is either a FH-representation 
$f(e^{i\th},\{\al_k\},\{\wt\bt_j\})$ 
of the symbol with $|||\wt\bt|||<1$, and (1.12) in \cite{DIKt1} applies, 
or there 
are at least 2 FH-representations with $|||\wt\bt|||=1$ (this is the situation for the example above).
In the latter case, the asymptotics
are given by Theorems 1.13 and 1.18 in \cite{DIKt1}, and the leading term is obtained as
a {\it sum} of contributions from different FH-representations of the symbol. 
Thus it can happen that these contributions sum up to zero. This is exactly the mechanism
by which the eigenvalues appear in the two examples below to which we restrict our attention from now on.
Of course, our considerations below can be easily generalized, but we feel it is more
useful to present the simplest cases which elucidate the mechanism. 
In the examples we consider, $f(z;\lb)$ has 2 Fisher-Hartwig singularities such that
$|||\bt|||=1$.
In the first example, the locations of the singularities of $f(z;\lb)$ will
depend on $\lb$; while in the second example, the locations are fixed, however, the
(imaginary parts of) $\bt$-parameters depend on $\lb$.

Before we describe our examples, we first recall some general facts about the spectra of 
Toeplitz operators and matrices. Let $f(e^{i\th})$ be a bounded, real-valued symbol on 
the unit circle, $f\in L^{\infty}(C)$. Let $T_n(f)$ be the associated Toeplitz matrix as above. 
Let $M_f$ denote the operator of multiplication by $f$ in $L^2(C)$ and let $T(f)$ be the 
Toeplitz operator associated with $f$ and acting in $\ell^2(0,1,2,\dots)$. All three operators
are self-adjoint and hence have real spectrum. The spectrum $\si(M_f)$ of $M_f$ is given by
the {\it essential range} of $f$,
\be\la{211}
\si(M_f)=\mbox{ ess range of } f=
\{\lb: \mbox{meas}\{\th: |f(e^{i\th})-\lb|<\ep\}>0 \mbox{ for all }\ep>0\}.
\ee
By a standard min-max argument,
\be\la{212}
\si(T_n(f))\subset [\inf\{\lb: \lb\in\si(T(f))\},\sup\{\lb: \lb\in\si(T(f))\}]
\ee
and
\be\la{213}
\si(T(f))\subset [\inf\{\lb: \lb\in\si(M_f)\},\sup\{\lb: \lb\in\si(M_f)\}]
\ee
and so by (\ref{211}),
\be\la{214}
\si(T(f))\subset [\essinf f, \esssup f].
\ee
By a theorem of Hartman and Wintner (see, e.g., \ci{BSbook}), we have equality in (\ref{214}):
\be\la{215}
\si(T(f))=[\essinf f, \esssup f].
\ee

As $T_n(f)$ converges strongly to $T(f)$ in $\ell^2(0,1,2,\dots)$ as $n\to\infty$, it follows 
by general principles (see, e.g., \ci{RS}) that all points in $\si(T(f))$ are limit points
of the spectra $\si(T_n(f))$, $n=1,2,\dots$, i.e. if $\lb\in\si(T(f))$, then 
$\lb=\lim_k \lb_{n_k}$, where $\lb_{n_k}\in \si(T_{n_k}(f))$. By the above considerations, we conclude
that
\be\la{216}
\overline{\cup_n\si(T_n(f))}=\si(T(f))=[\essinf f, \esssup f].
\ee
As the eigenvalues of $T_{n+1}(f)$ interlace with the eigenvalues of $T_n(f)$, the spectra $\si(T_n(f))$
fill in $\si(T(f))$ by casting, as it were, a finer and finer net.

Note that for each $n$, the eigenvalues $\lb^{(n)}_k$ of $T_n(f)$ lie in the open interval
$(\essinf f, \esssup f)$, apart from the trivial case $f(e^{i\th})\equiv \mbox{const}$ a.e.
Indeed, suppose that $(T_n(f)h)_j=\sum_{k=0}^{n-1}f_{j-k}h_k=s h_j$, $0\le j\le n-1$, where
$s=\esssup f$. Extending $h_j=0$ for $j<0$ and $j>n-1$, we have
\[
s \sum_{j=0}^{n-1}|h_j|^2=\sum_{j=0}^{n-1} \overline{h_j}\sum_{k=0}^{n-1}f_{j-k}h_k=
 \sum_{-\infty}^{\infty} \overline{h_j}\sum_{-\infty}^{\infty}f_{j-k}h_k=
\int_0^{2\pi}|h(\th)|^2 f(e^{i\th})\frac{d\th}{2\pi},
\]
where $h(\th)=\sum_{-\infty}^{\infty}e^{ij\th}h_j$. It follows that
\[
\int_0^{2\pi}(s-f(e^{i\th}))|h(\th)|^2\frac{d\th}{2\pi}=0.
\]
If $f$ is not identically a constant, then $s-f(e^{i\th})>0$ on a set of positive measure, and so
$h(\th)=\sum_{-\infty}^{\infty}e^{ij\th}h_j=0$ on a set of positive measure. 
Hence $h(\th)\equiv 0$ and therefore
$h_j=0$, $j=0,\dots,n-1$. A similar argument shows that $\lb^{(n)}_k>\essinf f$.

For our first example, we take $f$ as in (\ref{f1}). (Note that the $C^{\infty}$ condition on $V(z)$ 
can be relaxed, see \cite{DIKt2}). Theorem \ref{thm1} below was proved 
by B\"ottcher, Grudsky, and Maksimenko in \ci{BGM}
in the case of $e^{V(z)}$ being a trigonometric polynomial, using other methods. 
We have 
\begin{theorem}\la{thm1}
Let $f(z;0)\equiv f(z)$ be as described in and following (\ref{f1}). 
Assume furthermore that the second derivatives 
$f''(1)\neq 0$, $f''(e^{i\wt\th})\neq 0$,
and let $z_1=e^{i\th_1}$ and $z_2=e^{i\th_2}$, $0<\th_1<\th_2<2\pi$, be the zeros
of $f(z;\lb)=f(z;0)-\lb$,   $L <\lb< M$, where
$L=\min_{z\in C} f(z)$, $M=\max_{z\in C} f(z)$.

Then as $n\to\infty$, the eigenvalues $\lb^{(n)}_j$ of $T_n(f)$ satisfy
$L <\lb^{(n)}_j< M$ and
\be\la{cond0}
(n+1)\Psi(\lb^{(n)}_j)+\Theta(\lb^{(n)}_j)=\pi j+o(1),\qquad j=1,2,\dots,n,
\ee 
where 
\begin{multline}\la{PTf}
\Psi(\lb)=\Im(\ln f(z;\lb))_0={1\over 2}(\th_1-\th_2)+\pi,\\
\Theta(\lb)+\Psi(\lb)-{\pi\over 2}=
\Im \sum_{k=1}^\infty k (\ln f(z;\lb))_k (\ln f(z;\lb))_{-k}
=\Im \sum_{k=1}^\infty (z_1^k-z_2^k)(\ln|f(z;\lb)|)_k,
\end{multline}
with $(\cdot)_k$ denoting the $k$'th Fourier coefficient. 
The $o(1)$ term in (\ref{cond0}) is uniform in $j=1,2,\dots,n$. 
\end{theorem}

Relation (\ref{cond0}) leads to the following estimates on the eigenvalues $\lb^{(n)}_j$ of $T_n(f)$.

\begin{corollary}\la{cor1}
Let $f(z)$ be given as in Theorem \ref{thm1} and let $\lb^{(n)}_1<\lb^{(n)}_2<\cdots 
<\lb^{(n)}_n$ denote the eigenvalues of $T_n(f)$. Let $a_{\min}$, $a_{\max}$ be given 
as in (\ref{5114}) below, and let 
$0<\ep<a_{\min}/2$. 

\noindent (i) If $2\ep<j/n< 1-2\ep$, then for suitable constants $L<\lb_{\ep}<\mu_{\ep}<M$, 
as $n\to\infty$, $\lb_{\ep}<\lb^{(n)}_j<\mu_{\ep}$ and
\be\la{119}
\frac{c_1(\ep)}{n}\le  \lb^{(n)}_{j+1}-\lb^{(n)}_j\le \frac{c_2(\ep)}{n},
\ee
where the constants $0<c_1(\ep)<c_2(\ep)$ are uniform for $2\ep<j/n< 1-2\ep$.

\noindent (ii) If $0<j/n\le 2\ep$, then as $n\to\infty$,
\be\la{120}
\frac{j^2}{n^2a_{\max}^2}(1+o(1))\le \frac{\lb^{(n)}_j-L}{M-L}\le
\frac{\pi^2 j^2}{4n^2a_{\min}^2}(1+o(1))
\ee
and
\be\la{121}
\frac{c_3(\ep)}{n}\left(\frac{j}{n}\right)^2\le  
\lb^{(n)}_{j+1}-\lb^{(n)}_j\le \frac{c_4(\ep)}{n}\left(\frac{j}{n}\right)^2,
\ee
where the terms $o(1)$ and the constants $0<c_3(\ep)<c_4(\ep)$ are uniform for $0<j/n\le 2\ep$.
There are similar estimates for $j/n\ge 1-2\ep$, which correspond to replacing $j$ with $n-j$
in the above estimates.
\end{corollary}

We now turn to our second example.
Consider a real-valued function $f(z)$ such that there is a gap, say $(a,b)$, 
between components of the essential range of $f(z)$, $z\in C$.
By the preceding discussion, $[a,b]\subset\si(T(f))$ and each $\lb\in[a,b]$ is a limit
point of eigenvalues of the Toeplitz matrices $T_n(f)$.  However, as $(a,b)\cap\si(M_f)=\emptyset$,
we anticipate that these eigenvalues are sparsely distributed, and indeed, by \ci{Basor} (see also below),
for any subinterval $(a+\ep,b-\ep)$, $\ep>0$, the distance between the eigenvalues 
of $T_n(f)$ is of order $1/\ln n$ for $n$ sufficiently large.

In \cite{LS} (see also \cite{Sh})
Levitin and Shargorodsky considered a symbol of the form
$f_{\al}(e^{i\th})=-3/2+(\al/2)\cos(\sqrt{5}\th)$ for $\th\in[-\pi,s)$, 
$f_{\al}(e^{i\th})=2+\al\cos(\sqrt{2}\th)$ for $\th\in[s,\pi)$, where $s\in[0,\pi)$, $|\al|\le 1$,
and observed numerically (see, in particular, Figure 13 in  \cite{LS} and also Figure 1 in  \cite{LSS})
the following phenomenon of near-periodicity of the eigenvalues inside the gap of the range of $f_{\al}$.
For coprime integers $\ell,m\in\bbz\setminus\{0\}$, define
\be\la{321}
\om=\om(\ell,m)=
\begin{cases}
m, & \mbox{if $\ell$ and $m$ are odd}\cr
2m, & \mbox{if either $\ell$ or $m$ is even} 
\end{cases},
\ee
and if $\ell=0$, $m\in\bbz\setminus\{0\}$, let $\om=\om(0,m)=2$.
If $s=\pi \ell/m$, then for each eigenvalue
$\lb^{(n)}_k$ of $T_n(f_{\al})$ inside the interval $(-1,1)$, there appears to exist an eigenvalue
$\lb^{(n+\om)}_j$ of $T_{n+\om}(f_{\al})$ such that $|\lb^{(n)}_k-\lb^{(n+\om)}_j|=o(1/\ln (n+\om))$, i.e.,
$\lb^{(n+\om)}_j$ approaches $\lb^{(n)}_j$ faster than the logarithmic filling rate of the gap.

In \cite{LSS}, Levitin, Sobolev, and Sobolev considered the (modified) symbol
$f(e^{i\th})=-1$ for $\th\in[-\pi,s)$, $f(e^{i\th})=1$ for $\th\in[s,\pi)$, where $s\in[0,\pi)$,
and proved the near-periodicity of the eigenvalues of the square $T(f)^2$ in the gap $(0,1)$
when $s$ is a rational multiple of $\pi$ as above.  

Here we prove the near-periodicity conjecture of \cite{LS,LSS} assuming for simplicity that 
the range of the symbol $f(z)$ is 2 different real constants (such a symbol is of type (\ref{fFH})
with 2 jump-type singularities: see below). The proof, however, can be extended
to more general situations of type (\ref{fFH}), including the symbol $f_{\al}$ in \ci{LS} 
described above. By the above discussion, most of the eigenvalues of
$T_n(f)$ will be close to these constants, but we are interested in (the order $\ln n$)
eigenvalues which are in subintervals inside the gap.
We have 
  
\begin{theorem}\la{thm2}
Let $0\le \th_1<\th_2<2\pi$, $\ga>0$, and
\be\la{fLS}
f(e^{i\th})=
\begin{cases}
1,& \th\in[\th_2,2\pi)\cup[0,\th_1)\cr
e^{2\pi\ga},& \th\in[\th_1,\th_2)
\end{cases}.
\ee
Let, furthermore, $\th_1$, $\th_2$ be such that
\be\la{rational}
\th_2-\th_1=2\pi{p\over q},\qquad p,q\in\bbz,\qquad 0<p<q.
\ee
Consider the interval $I_{\ep}=(1+\ep,e^{2\pi\ga}-\ep)$ for a fixed $\ep$,
$0<\ep<(1+e^{2\pi\ga})/2$, and let $n$ be
sufficiently large. Then there exist constants $c_{\ell}>0$, $\ell=0,1,2$ which only depend on
$\ep$ and $\ga$, such that the distance between any 2 consecutive eigenvalues of $T_n(f)$
inside $I_{\ep}$ is bounded from below by $c_0/\ln n$, and from above, by $c_1/\ln n$. 
Any subinterval of $I_{\ep}$ of length $c_1/\ln n$ contains an eigenvalue.
For any eigenvalue $\lb^{(n)}_k$ of $T_n(f)$ inside $I_{\ep}$ there exists
an eigenvalue $\lb^{(n+q)}_j$ of $T_{n+q}(f)$ such that
\be\la{lblog}
|\lb^{(n)}_k-\lb^{(n+q)}_j|\le \frac{c_2}{n\ln n}.
\ee
\end{theorem}

\begin{remark}
Note that (\ref{rational}) encodes (\ref{321}) in a more compact way. Indeed, $\ell=0$ in (\ref{321})
corresponds to $\th_2-\th_1=\pi=2\pi{1\over 2}$. Thus $\om=2=q$. The case $\ell$, $m$ odd corresponds to
$\th_2-\th_1=\pi{m\pm\ell\over m}=2\pi{p\over m}$, as $m\pm\ell$ is even.
Thus $\om=m=q$. Finally, the case where either $\ell$ or $m$ is even corresponds to
$\th_2-\th_1=2\pi{m\pm\ell\over 2m}$, as  $m\pm\ell$ is odd. Thus $\om=2m=q$.
\end{remark}

\begin{remark}
Suppose $f(e^{i\th})=a$ for $\th\in[\th_1,\th_2)$, and $f(e^{i\th})=b$ for all other $\th\in[0,2\pi)$,
$b\neq a$. Then $\wt f(e^{i\th})\equiv (f(e^{i\th})-a)/(b-a)=0$ for $\th\in[\th_1,\th_2)$,
and $\wt f(e^{i\th})=1$ for all other $\th\in[0,2\pi)$.
As $T_n(\wt f)=(T_n(f)-a)/(b-a)$, it follows that the eigenvalues $\wt\lb$ of $T_n(\wt f)$ and
the eigenvalues $\lb$ of $T_n(f)$ are related through the elementary formula
$\wt\lb=(\lb-a)/(b-a)$. This clearly implies that the phenomenon of near-periodicity depends only
on $\th_2-\th_1$ and not on $a$ and $b$. It follows, in particular, that if the roles of
$[\th_1,\th_2)$ and its complement in $[0,2\pi)$ are reversed, the near-period $q$ should be the same.
As $2\pi-(\th_2-\th_1)=2\pi{q-p\over q}$, we see that this is indeed the case.
\end{remark}
At the end of the paper we discuss the relation of our results in Theorem \ref{thm2}
to a conjecture of Slepian \ci{Slepian} and its resolution in \ci{LW} by Landau and Widom.

\section{Asymptotics of some Toeplitz determinants}
Introduce the canonical Wiener-Hopf factorization of $e^{V(z)}$ (we assume $V(z)$ to be  
sufficiently smooth on $C$: see \cite{DIKt1} for details):
\begin{multline}\la{WienH}
e^{V(z)}=b_+(z) e^{V_0} b_-(z),\qquad b_+(z)=e^{\sum_{k=1}^\infty V_k z^k},
\qquad b_-(z)=e^{\sum_{k=-\infty}^{-1} V_k z^k},\\
V_k={1\over 2\pi}\int_0^{2\pi}V(e^{i\th})e^{-ik\th}d\th. 
\end{multline}
In the proofs of Theorem \ref{thm1} and \ref{thm2}, we will use the following formulae 
and asymptotic estimates for 
Toeplitz determinants, which are part of Theorems 1.1, 1.18, 1.8 in \cite{DIKt1}, and Theorem 1.1.
in \ci{DIKt2}.

\begin{lemma}\la{lemma1}
Let $F(z)$ be of the form (\ref{fFH}) with  $\Re\al_j>-{1\over 2}$, $\Re\bt_j\in (-1/2,1/2]$, 
$j=0,1,\dots,m$. Let the symbol $F^-(z)$ be obtained from $F(z)$ 
by replacing one $\bt_{j_0}$ with $\bt_{j_0}-1$ for some fixed $0\le j_0\le m$.
Then for sufficiently large $n$ ($n>N$ with some $N>0$), there exists a unique monic polynomial  
$\widehat\Phi_n(z)= z^n+\cdots$ of degree $n$ such that
\[\la{or0}
\int_0^{2\pi}\widehat\Phi_n(z^{-1})z^j F(z)d\theta=0,\qquad
z=e^{i\theta},\qquad j=0,1,\dots,n-1,
\]
and 
\be\la{DFDF}
D_n(F^{-})=z_{j_0}^n {\widehat\Phi_n(0)}D_n(F).
\ee
As $n\to\infty$,
\be\la{ashatphi}
{\widehat\Phi_n(0)}=
\sum_{j=0}^m n^{2\bt_j-1} z_j^{-n} \nu_j^{-1}
{\Gamma(1+\al_j-\bt_j)\over \Gamma(\al_j+\bt_j)}
\frac{b_-(z_j)}{b_+(z_j)}+o(1),
\ee
where
\be\la{nu}
\nu_j=\exp\left\{-i\pi\left(\sum_{p=0}^{j-1}\al_p-
\sum_{p=j+1}^m\al_p\right)\right\}
\prod_{p\neq j}\left({z_j\over z_p}\right)^{\al_p}
|z_j-z_p|^{2\bt_p}.
\ee
Here $1/\Gamma(\al_j+\bt_j)\equiv 0$ if $\al_j+\bt_j=0$.
If $V(z)$ is $C^{\infty}$ on the unit circle then the error term 
in (\ref{ashatphi}) $o(1)=O((n^{2|||\bt|||-2}+n^{-1})n^{2\max_j\Re\bt_j-1})$.

Furthermore,
\begin{multline}\la{asD}
D_n(F)=\exp\left[nV_0+\sum_{k=1}^\infty k V_k V_{-k}\right]
\prod_{j=0}^m b_+(z_j)^{-\al_j+\bt_j}b_-(z_j)^{-\al_j-\bt_j}\\
\times
n^{\sum_{j=0}^m(\al_j^2-\bt_j^2)}\prod_{0\le j<k\le m}
|z_j-z_k|^{2(\bt_j\bt_k-\al_j\al_k)}\left({z_k\over z_j e^{i\pi}}
\right)^{\al_j\bt_k-\al_k\bt_j}\\
\times
\prod_{j=0}^m\frac{G(1+\al_j+\bt_j) G(1+\al_j-\bt_j)}{G(1+2\al_j)}
\left(1+o(1)\right),
\end{multline}
where
$G(x)$ is Barnes' $G$-function. The double product over $j<k$ is set to $1$
if $m=0$.

If $V(z)$ is $C^{\infty}$ on the unit circle and $|||\bt|||<1$ then the error term 
in (\ref{asD}) is $o(1)=O(n^{|||\bt|||-1})$.
The error terms in (\ref{ashatphi}) and (\ref{asD}) are uniform 
in all $\al_j$, $\bt_j$ (and $N$ is independent of  $\al_j$, $\bt_j$)
for $\bt_j$ in compact subsets of the strip $\Re\bt_j\in(-1/2,1/2]$ and
for $\al_j$ in compact
subsets of the half-plane $\Re\al_j>-1/2$. These error terms are also uniform in the $z_j$'s
provided these points are at a fixed distance from one another on the unit circle, and uniform
in $V(z)$ (and $N$ is independent of $V(z)$)
provided the
$V_k$'s are uniformly bounded in absolute value by 
the Fourier coefficients of a sufficiently smooth function.
\end{lemma}

\section{Proof of Theorem \ref{thm1}}
As advertised above, the zeros $z_1$, $z_2$ of $f(z;\lb)$ can be regarded as giving rise to
Fisher-Hartwig singularities with $\al_1=\al_2=1/2$, $\bt_1=-\bt_2=1/2$. Thus 
$f(z;\lb)$ is of type (\ref{fFH}):
\begin{align}\la{fBGM}
f(z;\lb)=& f(z)-\lb=e^{V(z)}|z-z_1||z-z_2|g_{z_1,1/2}(z)g_{z_2,-1/2}(z)\left(\frac{z_1}{z_2}
\right)^{-1/2}\\
=&
e^{V(z)}4\sin{\th-\th_1\over 2}\sin{\th-\th_2\over 2}\left(\frac{z_1}{z_2}\right)^{-1/2}.\notag
\end{align}
By elementary calculus
\be\la{511}
R(e^{i\th};\lb)\equiv 
-\frac{f(z;\lb)}{4\sin{\th-\th_1\over 2}\sin{\th-\th_2\over 2}},
\ee
$\th\neq\th_1=\th_1(\lb)$, $\th_2=\th_2(\lb)$, extends to a continuous, strictly positive function
on $[0,2\pi)\times [L,M]$,
\be\la{512}
R(e^{i\th};\lb)\ge c>0.
\ee
We specify $V(z)=V(e^{i\th};\lb)$ uniquely by defining
\be\la{513}
V(z)=\ln R(e^{i\th};\lb)+{i\over 2}(\th_1-\th_2)+i\pi
\ee
for $(\th,\lb)\in [0,2\pi)\times [L,M]$, where $\ln$ denotes the principal branch.
Again by elementary calculus, using (\ref{512}), one sees that $V(e^{i\th};\lb)$ is a smooth function of
$\th$, with the property that each derivative $(\partial^{\ell}/\partial\th^{\ell})V(e^{i\th};\lb)$,
$\ell\ge 1$, is bounded uniformly for all $\lb\in [L,M]$:
\be\la{514}
\sup_{(\th,\lb)\in  [0,2\pi)\times [L,M]}\left|
{\partial^{\ell}\over\partial\th^{\ell}}V(e^{i\th};\lb)\right|\le c_{\ell}\qquad \ell=0,1,\dots
\ee

Now observe that $f(z;\lb)$ in (\ref{fBGM}) can be written as
\[
f(z;\lb)=F^{-}(z),
\]
where $F^{-}(z)$ satisfies the conditions of Lemma \ref{lemma1}, where
$j_0=2$, $m=2$, $\al_0=\bt_0=0$, $\al_1=\al_2=1/2$, $\bt_1=\bt_2=1/2$.

Note that the main asymptotic term (\ref{asD}) of $D_n(F)$  is non-zero
for $n$ sufficiently large, uniformly for $\lb\in [L+\ep,M-\ep]$, $\ep>0$, and the
condition for the eigenvalues of $T_n(f)$ (equivalently, the zeros of $D_n(F^{-})$)
comes from the vanishing of $\widehat\Phi_n(0)$ in (\ref{ashatphi}).
Thus, in our case the eigenvalues of $T_n(f)$ satisfy
\be\la{36}
\widehat\Phi_n(0)=
z_1^{-n} \nu_1^{-1}\frac{b_-(z_1)}{b_+(z_1)}+
z_2^{-n} \nu_2^{-1}\frac{b_-(z_2)}{b_+(z_2)}+o(1)=0,
\ee
where
\[
\nu_1=e^{i\pi/2}\left({z_1\over z_2}\right)^{1/2}|z_1-z_2|,\qquad
\nu_2=e^{-i\pi/2}\left({z_2\over z_1}\right)^{1/2}|z_1-z_2|,
\]
where the $o(1)$ term is uniform for $\lb\in[L+\ep,M-\ep]$ (cf. \ci{BGM}). 
Now as $\widehat\Phi_n(0)=\widehat\Phi_n(0;\lb)$ is a complex-valued function of $\lb\in[L+\ep,M-\ep]$,
and the eigenvalues of $T_n(f)$ are real, 
it is natural to consider a real-valued equivalent of (\ref{36}).
To this end we proceed as follows. As $f(e^{i\th};\lb)=f(e^{i\th})-\lb$ is real-valued
for $\lb$ real, we see from (\ref{DFDF}) that $z_2^n \widehat\Phi_n(0)D_n(F)=D_n(f-\lb)$
is real for $\lb\in[L+\ep,M-\ep]$. From (\ref{asD}), as $n\to\infty$,
\be
D_n(F)=\exp\left[nV_0+\sum_{k=1}^\infty k V_k V_{-k}\right](b_-(z_1)b_-(z_2))^{-1}(1+o(1)),
\ee
where the $o(1)$ term is uniform for $\lb\in[L+\ep,M-\ep]$. It follows, using 
$V_k=\overline{V_{-k}}$, $k\neq 0$, that
\begin{align*}
z_2^n \widehat\Phi_n(0)D_n(F)=&
i^{-1}\left[
\left({z_2\over z_1}\right)^{(n+1)/2}(b_-(z_2)b_+(z_1))^{-1}-
\left({z_1\over z_2}\right)^{(n+1)/2}(b_-(z_1)b_+(z_2))^{-1}+o(1)\right]\\
&\times |z_1-z_2|^{-1}\exp{
\left\{n\int_0^{2\pi}\ln R \frac{d\th}{2\pi}+\sum_{k=1}^\infty k |V_k|^2\right\}}
(1+o(1)),
\end{align*}
which can be written after an elementary calculation, and combining the $o(1)$ terms, as
\begin{align}\la{516}
D_n(f-\lb)|z_1-z_2|=&
2\left[\sin((n+1)\Psi(\lb)+\Th(\lb))+e_n(\lb)\right]\\
&\times (-1)^n e^{Z(\lb)}
\exp{\left\{n\int_0^{2\pi}\ln R \frac{d\th}{2\pi}+\sum_{k=1}^\infty k |V_k|^2
\right\}},\notag
\end{align}
where $\Psi$, $\Th$, and $Z$ are real-valued,
\begin{align}
\Psi(\lb)=& {1\over 2}(\th_1-\th_2)+\pi,\la{517}\\
\Th(\lb)=& {1\over 2i}\sum_{k=1}^{\infty}\left[V_k(z_1^k-z_2^k)-V_{-k}(z_1^{-k}-z_2^{-k})\right],
\la{518}\\
Z(\lb)=& -{1\over 2}\sum_{k=1}^{\infty}\left[V_k(z_1^k+z_2^k)+V_{-k}(z_1^{-k}+z_2^{-k})\right],
\la{519}
\end{align}
and 
\be\la{5110}
e_n(\lb)\mbox{  is  } o(1),\mbox{  uniformly for  $\lb\in[L+\ep,M-\ep]$}.
\ee
A more detailed analysis (see \ci{CK}) shows that in fact (\ref{5110}) holds uniformly
for all $\lambda \in (L,M)$. We will use this enhanced version of (\ref{5110}) in
what follows. 

From (\ref{516}) we see, in particular, that for $n$ sufficiently large
\begin{align}\la{5111}
E_n=& E_n(\lb)\equiv \sin((n+1)\Psi(\lb)+\Th(\lb))+e_n(\lb)=\\
=& D_n(f-\lb)\frac{(-1)^n}{2}|z_1-z_2| e^{-Z(\lb)}
\exp{\left\{-n\int_0^{2\pi}\ln R \frac{d\th}{2\pi}-\sum_{k=1}^\infty k |V_k|^2
\right\}}\notag
\end{align}
is real-valued and continuous, and
\be\la{5112}
E_n(\lb)=0\quad \Leftrightarrow \quad \mbox{$\lb$ is an eigenvalue of $T_n(f)$}.
\ee

Note that $\Psi(\lb)$ is clearly a strictly increasing continuous map 
from $[L,M]$ onto $[0,\pi]$. Furthermore, $\Psi(\lb)$ is smooth for $\lb\in (L,M)$,
and using the non-degeneracy of $f(e^{i\th})$ at its maximum and minimum,
one easily shows that
\be\la{5113}
\frac{d\Psi}{d\lb}=\frac{a(\lb)}{((\lb-L)(M-\lb))^{1/2}},
\ee
where $a(\lb)$ is continuous and positive on $[L,M]$,
\be\la{5114}
a_{\max}\ge a(\lb)\ge a_{\min}>0.
\ee

On the other hand, by our earlier discussion of $V(\lb)$, it follows that for each $k$,
$V_k=V_k(\lb)$ is continuous for $\lb\in[L,M]$.  Furthermore, taking $\ell=2$ in (\ref{514}), we
obtain the bound $\sup_{\lb\in[L,M]}|V_k(\lb)|\le c/k^2$, $k\neq 0$, and so 
$\Th(\lb)$ in (\ref{518}) is continuous on $[L,M]$. Since $z_1-z_2\to 0$ as $\lb\to L$ or $\lb\to M$, it
follows that 
\be\la{5115}
\lim_{\lb\downarrow L}\Th(\lb)=\lim_{\lb\uparrow M}\Th(\lb)=0.
\ee
Using the properties of $R(e^{i\th};\lb)$ in (\ref{511}), and taking $\ell=3$ in (\ref{514}),
we see by direct differentiation in (\ref{518}) that $\Th(\lb)$ is differentiable in $(L,M)$, and
\be\la{5116}
\frac{d\Th}{d\lb}=\frac{b(\lb)}{((\lb-L)(M-\lb))^{1/2}},
\ee
where $b(\lb)$ is continuous on $[L,M]$,
\be\la{5117}
b_{\max}\ge b(\lb)\ge b_{\min}.
\ee

Hence
\be\la{5118}
G(\lb)\equiv
 (n+1)\Psi(\lb)+\Th(\lb)
\ee
is differentiable in $(L,M)$, and for $n$ sufficiently large
\[
\frac{dG}{d\lb}=\frac{(n+1)a(\lb)+b(\lb)}{((\lb-L)(M-\lb))^{1/2}}>0
\]
as $(n+1)a(\lb)+b(\lb)\ge(n+1)a_{\min}+b_{\min}$ and $a_{\min}>0$.
Thus for $n$ sufficiently large, $G(\lb)$ is a strictly increasing map from $[L,M]$ onto
$[G(L),G(M)]=[0,(n+1)\pi]$. Hence there exist unique points 
\[
L<\wh\lb^{(n)}_1<\wh\lb^{(n)}_2<\cdots<\wh\lb^{(n)}_n<M
\]
such that
\[
G(\wh\lb^{(n)}_j)=j\pi,\qquad 1\le j\le n.
\]
Set
\be\la{5119}
\ep_n=2\sup_{\lb\in [L,M]}|e_n(\lb)|,\qquad \wh\ep_n=\arcsin\ep_n>0.
\ee
By (\ref{5110}), $\ep_n$, and hence $\wh\ep_n$, converge to zero as $n\to\infty$. Again there exist 
unique points $\wh\lb^{(n)}_{j,\pm}$ such that 
\[
G(\wh\lb^{(n)}_{j,\pm})=j\pi\pm\wh\ep_n.
\]
For $n$ sufficiently large, by the monotonicity of $G(\lb)$,
\[
L<\wh\lb^{(n)}_{1,-}<\wh\lb^{(n)}_1<\wh\lb^{(n)}_{1,+}<
\wh\lb^{(n)}_{2,-}<\wh\lb^{(n)}_2<\wh\lb^{(n)}_{2,+}<
\cdots<\wh\lb^{(n)}_{n,-}<\wh\lb^{(n)}_n<\wh\lb^{(n)}_{n,+}<M
\]
Now
\[
E_n(\wh\lb^{(n)}_{j,\pm})=\sin(j\pi\pm\wh\ep_n)+e_n(\wh\lb^{(n)}_{j,\pm})=
\pm(-1)^j\ep_n+e_n(\wh\lb^{(n)}_{j,\pm}),
\]
and it follows from (\ref{5119}) that $E_n(\wh\lb^{(n)}_{j,+})$ and $E_n(\wh\lb^{(n)}_{j,-})$
have opposite signs, and hence by the continuity of (the real-valued function) $E_n(\lb)$, that
\be\la{5120}
E_n(\lb^{(n)}_j)=0,\quad\mbox{for some points}\quad  
  \wh\lb^{(n)}_{j,-}<\lb^{(n)}_j<\wh\lb^{(n)}_{j,+},\quad j=1,\dots,n.
\ee
Again by the monotonicity of $G(\lb)$,
\[
j\pi-\wh\ep_n=G(\wh\lb^{(n)}_{j,-})<G(\lb^{(n)}_j)<G(\wh\lb^{(n)}_{j,+})=j\pi+\wh\ep_n,
\]
i.e.,
\be\la{5121}
(n+1)\Psi(\lb^{(n)}_j)+\Th(\lb^{(n)}_j)=j\pi+o(1),\qquad j=1,\dots,n,
\ee
where the error term $o(1)$ is uniform in $j$. As $T_n(f)$ has precisely $n$ eigenvalues,
it follows from the pigeon hole principle that the points $\lb^{(n)}_j\in 
(\wh\lb^{(n)}_{j,-},\wh\lb^{(n)}_{j,+})$ with $E_n(\lb^{(n)}_j)=0$ are unique and comprise all the
eigenvalues of $T_n(f)$. This completes the proof of Theorem \ref{thm1} up to the formulae in
(\ref{PTf}) which relate $\Psi$ and $\Th$ directly to $f$. However, from (\ref{511}) and (\ref{513})
\be
V_k=(\ln|f(z;\lb)|)_k-(\ln|z-z_1||z-z_2|)_k,\qquad k\neq 0,
\ee
and by an elementary calculation,
\be
(\ln|z-z_j|)_k=-{1\over 2|k|}z_j^{-k},\qquad j=1,2,\qquad k\neq 0.
\ee
Substitution into (\ref{518}) yields for all $\lb\in (L,M)$
\be\la{Theta}
\Th(\lb)=\Im\left(\sum_{k=1}^{\infty}(z_1^k-z_2^k)
(\ln|f(z;\lb)|)_k\right)-\frac{1}{2}(\pi+\th_1-\th_2)=
\Im \sum_{k=1}^\infty (z_1^k-z_2^k)(\ln|f(z;\lb)|)_k-\Psi(\lb)+{\pi\over 2},
\ee
which establishes one of the expressions in (\ref{PTf}).
The final equation in (\ref{PTf}),
\[
\Im \sum_{k=1}^\infty k (\ln f(z;\lb))_k (\ln f(z;\lb))_{-k}
=\Im \sum_{k=1}^\infty (z_1^k-z_2^k)(\ln|f(z;\lb)|)_k,
\]
is straightforward to verify. Here, by (\ref{fBGM}), (\ref{513}),
\be\la{lnf}
\ln f(e^{i\th};\lb)=
\begin{cases}
\ln|f(z;\lb)|+i\pi,& \th\in(\th_2,2\pi)\cup[0,\th_1)\cr
\ln|f(z;\lb)|,& \th\in(\th_1,\th_2)
\end{cases}.
\ee
$\Box$

\section{Proof of Corollary \ref{cor1}}
Corollary \ref{cor1} provides detailed information on the 
behavior of the eigenvalues $\lb^{(n)}_j$ as $n\to\infty$. We have
\be\la{5122}
\Psi(\lb^{(n)}_j)={j\over n}\pi+O\left({1\over n}\right),
\ee
where the error term $O(1/n)$ is uniform for $j=1,2,\dots,n$. Let $\ep>0$ be small and given.
Then there exist $L<\lb_{\ep}<\mu_{\ep}<M$ such that $\Psi(\lb_{\ep})=\pi\ep$,
$\Psi(\mu_{\ep})=\pi(1-\ep)$. Suppose $2\ep<j/n<1-2\ep$. Then by (\ref{5122}), as $n\to\infty$, 
$\pi\ep<\Psi(\lb^{(n)}_j)<\pi(1-\ep)$ and so by monotonicity $\lb_{\ep}<\lb^{(n)}_j<\mu_{\ep}$. 
To estimate the spacing $\lb^{(n)}_{j+1}-\lb^{(n)}_j$ for such $j$, write
\[
G(\lb^{(n)}_{j+1})-G(\lb^{(n)}_j)=\pi+o(1).
\]
On the other hand,
\be\la{5123}
G(\lb^{(n)}_{j+1})-G(\lb^{(n)}_j)=G'(\xi^{(n)}_j)(\lb^{(n)}_{j+1}-\lb^{(n)}_j)=
\frac{(n+1)a(\xi^{(n)}_j)+b(\xi^{(n)}_j)}{((\xi^{(n)}_j-L)(M-\xi^{(n)}_j))^{1/2}}
(\lb^{(n)}_{j+1}-\lb^{(n)}_j),
\ee
for some $\lb_{\ep}<\lb^{(n)}_{j}<\xi^{(n)}_j<\lb^{(n)}_{j+1}<\mu_{\ep}$. We have
\[
\frac{(n+1)a_{\min}+b_{\min}}{((\mu_{\ep}-L)(M-\lb_{\ep}))^{1/2}}\le
\frac{(n+1)a(\xi^{(n)}_j)+b(\xi^{(n)}_j)}{((\xi^{(n)}_j-L)(M-\xi^{(n)}_j))^{1/2}}\le
\frac{(n+1)a_{\max}+b_{\max}}{((\lb_{\ep}-L)(M-\mu_{\ep}))^{1/2}},
\]
which yields (\ref{119}) with the constants $0<c_1(\ep)<c_2(\ep)$ independent of $j$ 
for $2\ep<j/n<1-2\ep$ and sufficiently large $n$.

Now suppose $0<j/n\le 2\ep$; the case $1\ge j/n\ge 1-2\ep$ is similar.
We have
\[
G(\lb^{(n)}_j)=G(\lb^{(n)}_j)-G(L)=j\pi+o(1).
\]
On the other hand,
\[
G(\lb^{(n)}_j)-G(L)=\int_0^1 G'(\lb(s))(\lb^{(n)}_j-L)ds=
(\lb^{(n)}_j-L)\int_0^1\frac{(n+1)a(\lb(s))+b(\lb(s))}{[(\lb(s)-L)(M-\lb(s))]^{1/2}}ds,
\]
where
$\lb(s)=L+s(\lb^{(n)}_j-L)$, and so
\begin{align*}
((n+1)a_{\min}+b_{\min})(\lb^{(n)}_j-L)
\int_0^1\frac{ds}{[(\lb(s)-L)(M-\lb(s))]^{1/2}}\le
G(\lb^{(n)}_j)-G(L)\\
\le
((n+1)a_{\max}+b_{\max})(\lb^{(n)}_j-L)
\int_0^1\frac{ds}{[(\lb(s)-L)(M-\lb(s))]^{1/2}},
\end{align*}
and after integration we obtain
\be\la{5125}
\frac{\pi j+o(1)}{(n+1)a_{\max}+b_{\max}}\le
\arcsin\frac{2\lb^{(n)}_j-(L+M)}{M-L}+{\pi\over 2}\le
\frac{\pi j+o(1)}{(n+1)a_{\min}+b_{\min}}.
\ee
We may assume without loss that $\ep<a_{\min}/2$. Then for $n$ sufficiently large,
\[
\frac{\pi j+o(1)}{(n+1)a_{\min}+b_{\min}}<\frac{\pi 2\ep}{a_{\min}}<\pi,
\]
and we obtain from (\ref{5125})
\[
\cos \frac{\pi j+o(1)}{(n+1)a_{\max}+b_{\max}}\ge 
-\frac{2\lb^{(n)}_j-(L+M)}{M-L}\ge
\cos \frac{\pi j+o(1)}{(n+1)a_{\min}+b_{\min}}.
\]
Therefore
\be
\sin^2 \frac{\pi j+o(1)}{2((n+1)a_{\max}+b_{\max})}\le
\frac{\lb^{(n)}_j-L}{M-L}\le
\sin^2 \frac{\pi j+o(1)}{2((n+1)a_{\min}+b_{\min})}.
\ee
Now for $\ep<a_{\min}/2$ as $n\to\infty$,
\[
0<\frac{\pi j+o(1)}{2((n+1)a_{\max}+b_{\max})}<
\frac{\pi j+o(1)}{2((n+1)a_{\min}+b_{\min})}<\frac{\pi}{2},
\]
and since $2/\pi\le(\sin x)/x\le 1$ for $0\le x\le\pi/2$, we obtain
(\ref{120}) uniformly for $j/n<2\ep<a_{\min}$, $n\to\infty$.
Recalling that 
$\lb^{(n)}_j<\xi^{(n)}_j<\lb^{(n)}_{j+1}$, and then utilizing (\ref{120}), we obtain 
(\ref{121}) from (\ref{5123}) uniformly for $j/n\le 2\ep<a_{\min}$, $n\to\infty$. $\Box$

\section{Proof of Theorem \ref{thm2}}
Let $\th_1>0$. Then
the function (\ref{fLS}) can be written in the form (\ref{fFH})
with  $m=2$, $z_1=e^{i\th_1}$, $z_2=e^{i\th_2}$,
$\al_0=\bt_0=0$, $\al_1=\al_2=0$,
\[
\bt_1=i\ga,\qquad \bt_2=-\bt_1,\qquad e^{V(z)}=\left({z_1\over z_2}\right)^{\bt_1}=e^{V_0},
\]
i.e.,
\be
f(z)=g_{z_1,i\ga}(z)g_{z_2,-i\ga}(z).
\ee
In the case of $\th_1=0$, $f(z)$ is written similarly:
we just need to replace indices $1$ by $0$, $2$ by $1$, and set $m=1$ in (\ref{fFH}).
For simplicity, we assume from now on that $\th_1>0$.

Note that $f(z;\lb)$, $\lb\in(1,e^{2\pi\ga})$, is also of type (\ref{fFH}) 
with the same points of singularities $z_1$, $z_2$ but 
with $\bt$ parameters and $V(z)$ now depending on $\lb$, namely
\be
f(z;\lb)=f(z)-\lb= e^{V^{(\lb)}(z)}g_{z_1,\bt_1^{(\lb)}}(z)g_{z_2,-\bt_1^{(\lb)}}(z)
\left({z_1\over z_2}\right)^{-\bt_1^{(\lb)}},\qquad \lb\in(1,e^{2\pi\ga}),
\ee
where
\be
\bt_1^{(\lb)}=i\ga^{(\lb)}+{1\over 2},\qquad e^{2\pi\ga^{(\lb)}}=
\frac{e^{2\pi\ga}-\lb}{\lb-1},
\ee
and
\be
e^{V^{(\lb)}(z)}=e^{i\pi}(\lb-1)\left({z_1\over z_2}\right)^{\bt_1^{(\lb)}}=e^{V^{(\lb)}_0}.
\ee

Note now that this
\[
f(z;\lb)=F^{-}(z),
\]
where $F^{-}(z)$ satisfies the conditions of Lemma \ref{lemma1}, where $j_0=2$, $m=2$, 
$\al_0=\bt_0=0$, $\al_1=\al_2=0$, $\bt_1=i\ga^{(\lb)}+1/2$, $\bt_2=-i\ga^{(\lb)}+1/2$.
From (\ref{asD}), as $n\to\infty$,
\be
\begin{aligned}
D_n(F)=&(1-\lb)^n|z_1-z_2|^{2(\ga^{(\lb)})^2+1/2}\left|G\left({1\over 2}+i\ga^{(\lb)}\right)\right|^2
\left|G\left({3\over 2}+i\ga^{(\lb)}\right)\right|^2
n^{2(\ga^{(\lb)})^2-1/2}\\
&\times \left({z_1\over z_2}\right)^{n/2}e^{-(\th_1-\th_2)n\ga^{(\lb)}}
\left[1+O\left({1\over n}\right)\right]\neq 0,
\end{aligned}
\ee
and the
condition for the eigenvalues of $T_n(f)$ (equivalently, for $D_n(F^{-})=0$)
comes from the vanishing of $\widehat\Phi_n(0)$ in (\ref{ashatphi}):
\be
z_1^{-n} (n|z_1-z_2|)^{2i\ga^{(\lb)}}
\frac{\Gamma(1/2-i\ga^{(\lb)})}{\Gamma(1/2+i\ga^{(\lb)})}+
z_2^{-n} (n|z_1-z_2|)^{-2i\ga^{(\lb)}}
\frac{\Gamma(1/2+i\ga^{(\lb)})}{\Gamma(1/2-i\ga^{(\lb)})}+O\left({1\over n}\right)=0,
\ee
where the $O(n^{-1})$ term is uniform for $\lb\in I_{\ep/2}$. 
As in Theorem \ref{thm1}, however, $\wh\Phi_n(0)$ is not real-valued and we again consider instead
the real-valued combination $z_2^n\wh\Phi_n(0)D_n(F)=D_n(f-\lb)$. Using the above asymptotics, and 
combining the $O(n^{-1})$ error terms, we obtain as $n\to\infty$
\[
z_2^n\wh\Phi_n(0)D_n(F)=P_n(\lb)E_n(\lb),
\]
where $P_n(\lb)$ is real-valued and non-zero and
\be
\begin{aligned}
E_n(\lb)=&{1\over 2}\left(\frac{z_2}{z_1}\right)^{n/2} (n|z_1-z_2|)^{2i\ga^{(\lb)}}
\frac{\Gamma(1/2-i\ga^{(\lb)})}{\Gamma(1/2+i\ga^{(\lb)})}\\
&+
{1\over 2}\left(\frac{z_1}{z_2}\right)^{n/2}
 (n|z_1-z_2|)^{-2i\ga^{(\lb)}}
\frac{\Gamma(1/2+i\ga^{(\lb)})}{\Gamma(1/2-i\ga^{(\lb)})}+e_n(\lb),
\end{aligned}
\ee
where $e_n(\lb)=O(n^{-1})$ uniformly for $\lb\in I_{\ep/2}$.
Set
\be
h(\lb)=\arg\Gamma(1/2+i\ga^{(\lb)}).
\ee
As the r.h.s. here is smooth and nonzero, $h(\lb)$ is uniquely determined as a smooth
function on $(1,e^{2\pi\ga})$ with $h({1\over 2}(1+e^{2\pi\ga}))= \arg\Gamma(1/2)=0$.
In terms of $h(\lb)$,
\be
E_n(\lb)=\cos\left({n\over 2}(\th_2-\th_1)+H_n(\lb)\right)+e_n(\lb),
\ee
where
\be
H_n(\lb)=2\ga^{(\lb)}\ln(n|z_1-z_2|)-2 h(\lb).
\ee
We have that $E_n(\lb)$ is real-valued and for $\lb\in I_{\ep}\subset  I_{\ep/2}$,
\be
E_n(\lb)=0\quad \Leftrightarrow \quad \mbox{$\lb$ is an eigenvalue of $T_n(f)$}.
\ee

Now the derivatives  ${d\over d\lb}\ga^{(\lb)}$,
${d\over d\lb}h(\lb)$ are bounded on  $I_{\ep/2}$,
\be\la{811}
c_{\min}\le {d\over d\lb}\ga^{(\lb)} \le c_{\max}<0,\qquad
d_{\min}\le {d\over d\lb}h(\lb) \le d_{\max},
\ee
and hence, for $n$ sufficiently large, $H_n(\lb)$ is strictly monotone, and maps
$I_{\ep}=[1+\ep, e^{2\pi\ga}-\ep]$ bijectively onto 
$[H_n(e^{2\pi\ga}-\ep),H_n(1+\ep)]$ which is of length
\[
2(\ga^{(1+\ep)}-\ga^{(e^{2\pi\ga}-\ep)})\ln(n|z_1-z_2|)-2 (h(1+\ep)-h(e^{2\pi\ga}-\ep)).
\]
In particular, there are $O(\ln n)$ values of $k\in\bbz$ such that
\be
k+{1\over 2}\in \left[\frac{\th_2-\th_1}{2\pi}n+{1\over\pi}H_n(e^{2\pi\ga}-\ep), 
\frac{\th_2-\th_1}{2\pi}n+{1\over\pi}H_n(1+\ep)
\right].
\ee

Let
\be
K=K_n=\{k_{\min}(n),k_{\min}(n)+1,\dots,k_{\max}(n)\},
\ee
where for each $k\in K$ there exists a unique point $\wh\lb^{(n)}_k\in [1+\ep,e^{2\pi\ga}-\ep]$
such that $\frac{\th_2-\th_1}{2\pi}n+{1\over\pi}H_n(\wh\lb^{(n)}_k)=k+{1\over 2}$. 
If $\wh\lb^{(n)}_{k_{\min}}<\lb<e^{2\pi\ga}-\ep/2$, then for some 
$\wh\lb^{(n)}_{k_{\min}}<\xi_n<\lb$,
\[
H_n(\wh\lb^{(n)}_{k_{\min}})-H_n(\lb)=
H_n'(\xi_n)(\wh\lb^{(n)}_{k_{\min}}-\lb)\ge
c_5(\lb- \wh\lb^{(n)}_{k_{\min}})\ln n,
\]
for some $c_5>0$, $n$ sufficiently large. Choosing $n$ even larger, if necessary, to ensure that
$c_5{\ep\over 2}\ln n >1$, and recalling that 
$\wh\lb^{(n)}_{k_{\min}}\le e^{2\pi\ga}-\ep<e^{2\pi\ga}-\ep/2$, we conclude that there exists
$\wh\lb^{(n)}_{k_{\min}}<\wh\lb<e^{2\pi\ga}-\ep/2$ such that
\[
\frac{\th_2-\th_1}{2\pi}n+{1\over\pi}H_n(\wh\lb)=
\frac{\th_2-\th_1}{2\pi}n+{1\over\pi}H_n(\wh\lb^{(n)}_{k_{\min}})-1=(k_{\min}-1)+{1\over 2}.
\]
Write $\wh\lb=\wh\lb^{(n)}_{k_{\min}-1}$.

By a similar argument there exists $\wh\lb^{(n)}_{k_{\max}+1}$,
$1+\ep/2<\wh\lb^{(n)}_{k_{\max}+1}<\wh\lb^{(n)}_{k_{\max}}$ such that
\[
\frac{\th_2-\th_1}{2\pi}n+{1\over\pi}H_n(\wh\lb^{(n)}_{k_{\max}+1})=(k_{\max}+1)+{1\over 2}.
\]
These arguments show that for $n$ sufficiently large
$\{k_{\min}-1/2,k_{\min}+1/2,k_{\max}+1/2,k_{\max}+3/2\}$ all lie in
\[
\left[ 
\frac{\th_2-\th_1}{2\pi}n+{1\over\pi}H_n(e^{2\pi\ga}-{\ep\over 2}),
\frac{\th_2-\th_1}{2\pi}n+{1\over\pi}H_n(1+{\ep\over 2})\right].
\]
It now follows as in the proof of Theorem \ref{thm1} that there exist points
\begin{align*}
1+{\ep\over 2}<\wh\lb^{(n)}_{k_{\max}+1}<
\wh\lb^{(n)}_{k_{\max},-}<\wh\lb^{(n)}_{k_{\max}}<\wh\lb^{(n)}_{k_{\max},+}<
\wh\lb^{(n)}_{k_{\max}-1,-}<\wh\lb^{(n)}_{k_{\max}-1}<\wh\lb^{(n)}_{k_{\max}-1,+}<\cdots\\
<\wh\lb^{(n)}_{k_{\min},-}<\wh\lb^{(n)}_{k_{\min}}<\wh\lb^{(n)}_{k_{\min},+}<
\wh\lb^{(n)}_{k_{\min}-1}<e^{2\pi\ga}-{\ep\over 2}
\end{align*}
such that each of the intervals $(\wh\lb^{(n)}_{k,-},\wh\lb^{(n)}_{k,+})$, $k_{\min}\le k\le
k_{\max}$, contains a zero $\lb^{(n)}_k$ of $E_n$, i.e. an eigenvalue of $T_n(f)$. 

In contrast to the proof of Theorem \ref{thm1}, there is no convenient pigeon hole principle 
to apply, so that a priori there could be more than one zero $\lb^{(n)}_k$ of $E_n$ in the interval 
 $(\wh\lb^{(n)}_{k,-},\wh\lb^{(n)}_{k,+})$. However, the eigenvalues of $T_n(f)$ and $T_{n+1}(f)$ 
interlace and therefore if there were two or more eigenvalues of $T_n(f)$ in that interval, then
the interval would also have to contain an eigenvalue of $T_{n+1}(f)$, but proceeding as in the proof
of Theorem \ref{thm1}, we obtain $E_n(\wh\lb^{(n)}_{k,\pm})=O(n^{-1})$ uniformly for
$k\in K_n$, and hence by monotonicity, $E_n(\lb)=O(n^{-1})$ for any $\lb\in(\wh\lb^{(n)}_{k,-},
\wh\lb^{(n)}_{k,+})$, $k_{\min}\le k\le k_{\max}$. For such $\lb$, as $n\to\infty$,
\begin{align*}
E_{n+1}(\lb)=&\cos \left[
{n\over 2}(\th_2-\th_1)+ H_n(\lb)+{\pi p\over q} +O\left({1\over n}\right)\right]
+e_{n+1}(\lb)\\
=& E_n(\lb)\cos\left[{\pi p\over q} +O\left({1\over n}\right)\right]-
\sin \left[
\frac{n}{2}(\th_2-\th_1)+ H_n(\lb)\right]\sin\left[{\pi p\over q} +O\left({1\over n}\right)\right]
+O\left({1\over n}\right)\\
=&
\pm\sin{\pi p\over q}+O\left({1\over n}\right).
\end{align*}
However, as $0<p/q<1$, it follows that for $n$ sufficiently large, $E_{n+1}(\lb)$ has no zeros in
$(\wh\lb^{(n)}_{k,-},\wh\lb^{(n)}_{k,+})$, and thus these intervals contain one, and only one,
eigenvalue of $T_n(f)$. The labelling $\lb^{(n)}_k$ for the eigenvalues of  $T_n(f)$ by the integer
$k$ is therefore appropriate.

To estimate the gap between eigenvalues $\lb^{(n)}_k$, $\lb^{(n)}_{k+1}$ with $k,k+1\in K_n$,
we note that as in Theorem \ref{thm1},
\be\la{lbkn}
\frac{\th_2-\th_1}{2}n+ H_n(\lb^{(n)}_k)=(k+{1\over 2})\pi+O(n^{-1}),
\ee
and similarly,
$\frac{\th_2-\th_1}{2}n+ H_n(\lb^{(n)}_{k+1})=(k+{3\over 2})\pi+O(n^{-1})$, and so
\[
\pi+O\left({1\over n}\right)= H_n(\lb^{(n)}_{k+1})-H_n(\lb^{(n)}_k)=
H_n'(\xi^{(n)}_k)(\lb^{(n)}_{k+1}-\lb^{(n)}_k),
\]
where $\xi^{(n)}_k\in (\lb^{(n)}_{k+1},\lb^{(n)}_k)\in I_{\ep/2}$.

Using (\ref{811}) we now obtain the bounds
\be
{c_0\over\ln n}\le \lb^{(n)}_k-\lb^{(n)}_{k+1}\le {c_1\over\ln n}
\ee
for suitable constants $0<c_0=c_0(\ep,\ga)<c_1=c_1(\ep,\ga)$.

Finally, using (\ref{rational}) and (\ref{lbkn}), we note that for $k\in K_n$,
\be
\begin{aligned}
\frac{\th_2-\th_1}{2\pi}(n+q)+{1\over\pi}H_{n+q}(\lb^{(n)}_k)=&
\frac{\th_2-\th_1}{2\pi}n+{1\over\pi}H_n(\lb^{(n)}_k)+ p + O\left({1\over n}\right)\\
=&k+p+{1\over 2}+ O\left({1\over n}\right).
\end{aligned}
\ee
It now follows easily from our previous calculation that for $n$ sufficiently large, there exists 
an eigenvalue
$\lb^{(n+q)}_{k+p}\in I_{\ep/2}$ of $T_{n+q}(f)$ such that 
$|\lb^{(n+q)}_{k+p}-\lb^{(n)}_k|\le\frac{c_2}{n\ln n}$. This completes the proof of Theorem \ref{thm2}.
$\Box$

\section{On a conjecture of Slepian}
Let $T=(t_1,t_2)$ and $S=(s_1,s_2)$ be intervals in $\bbr$ and let $h_T(x)={1\over 2\pi}
\int_T e^{ix\om}d\om$. Let $A_{S,T}(c)$, $c>0$, denote the operator
\[
(A_{S,T}(c)f)(x)=\int_{cS}h_T(x-y)f(y)dy,\qquad x\in cS,
\]
acting on $L^2(cS)$. In \ci{LW}, Landau and Widom consider the asymptotics of the eigenvalues 
$\{\lb_k(c)\}$ of the operator $A_{S,T}(c)$ as $c\to\infty$. Using trace class methods they prove, 
in particular, the following conjecture of Slepian \ci{Slepian}. 
Suppose
\be\la{90}
k=\left[{1\over 2\pi} |S||T|c+{1\over\pi^2} b \ln c\right],\qquad b\in\bbr.
\ee
Then $\lb_k(c)\to (1+e^b)^{-1}$ as $c\to\infty$.

To relate this result to the results in this paper, let $T=(\th_1,\th_2)$ and $S=(0,s_2)$. Discretizing 
the eigenvalue equation $A_{S,T}(c)\phi_k=\lb_k(c)\phi_k$, we obtain
\be\la{91}
\sum_{j=0}^{n-1}h_T(\ell\de-j\de)\phi_k(j\de)\de\sim
\lb_k(c)\phi_k(\ell\de),
\ee
where $\de$ is small and positive and $0\le\ell\le n-1$ with $n=[cs_2/\de]$. 
Note that
$h_T(j\de)=\int_T e^{ij\de\om}{d\om\over 2\pi}={1\over\de}\int_{\de\th_1}^{\de\th_2}e^{ij\th}
{d\th\over 2\pi}$, and so (\ref{91}) corresponds to the eigenvalue problem for (the transpose of) 
the Toeplitz matrix $T_n(\chi_\de)$, where $\chi_\de$ is the characteristic function of
the interval $(\de\th_1,\de\th_2)$. The symbol $\chi_\de$ corresponds to $f$ in (\ref{fLS}) provided
we replace $\lb$ by $\lb-1$ and choose $\ga$ such that $e^{2\pi\ga}=2$.
Now from the proof of Theorem \ref{thm2},
\[
k+{1\over 2}=  \frac{\th_2-\th_1}{2\pi}n\de +{1\over\pi}
\left(2\ga^{(\lb^{(n)}_k)}\ln(n|z_1-z_2|)-2 h(\lb^{(n)}_k)\right)+O\left({1\over n}\right),
\]
i.e.,
\be\la{92}
k={1\over 2\pi} |S||T|c+{2\over\pi}\ga^{(\lb^{(n)}_k)} \ln c+O(1).
\ee
Subtracting (\ref{90}) from (\ref{92}) we see that $2\pi\ga^{(\lb^{(n)}_k)}-b=O((\ln c)^{-1})$.
Thus
\[
e^{2\pi\ga^{(\lb^{(n)}_k)}}=\frac{1-(\lb^{(n)}_k-1)}{\lb^{(n)}_k-1}\to e^b,
\]
i.e. $\lb^{(n)}_k-1\to (1+e^b)^{-1}$, which is Slepian's formula. To make these arguments rigorous,
we must control all estimates uniformly as $\de\downarrow 0$; this can probably be done using the methods
in \ci{DIKZ} where an analogous uniformity problem arises, but we provide no further details.

\section*{Acknowledgements}
P. Deift was supported
in part by NSF grant \# DMS 1001886.
A. Its was supported
in part by NSF grant \#DMS-1001777.
I. Krasovsky was
supported in part by EPSRC grant \#EP/E022928/1. 
We are grateful to Michael Levitin and Alexander Sobolev for stimulating our interest 
in the Toeplitz eigenvalue problem. 
We would also like to thank
Henry Landau for drawing our attention to \ci{LW} and for some very useful remarks.



\end{document}